\documentclass[12pt,english,final]{elsarticle}
\usepackage{mathpazo}

\usepackage[T1]{fontenc}
\usepackage[latin9]{inputenc}
\usepackage{geometry}
\usepackage{color}
\usepackage{enumitem}
\usepackage{amsmath}
\usepackage{amssymb}
\usepackage{cancel}
\usepackage{graphicx}
\PassOptionsToPackage{normalem}{ulem}
\usepackage{ulem}

\makeatletter

\@ifundefined{showcaptionsetup}{}{%
 \PassOptionsToPackage{caption=false}{subfig}}
\usepackage{subfig}
\makeatother

\usepackage{babel}
\begin{document}

\title{An asymptotic-preserving semi-Lagrangian algorithm for the anisotropic
heat transport equation with arbitrary magnetic fields}

\author[lanl]{L. Chacón\corref{cor1}}

\ead{chacon@lanl.gov}

\author[epfl]{G. Di Giannatale}

\cortext[cor1]{Corresponding author}

\address[lanl]{Los Alamos National Laboratory, Los Alamos, NM 87545, USA}

\address[epfl]{Ecole Polytechnique Fédérale de Lausanne, CH-1015 Lausanne, Switzerland}
\begin{abstract}
We extend the recently proposed semi-Lagrangian algorithm for the
extremely anisotropic heat transport equation {[}Chacón et al.,\emph{
J. Comput. Phys.}, \textbf{272} (2014){]} to deal with arbitrary magnetic
field topologies. The original scheme (which showed remarkable numerical
properties) was valid for the so-called tokamak-ordering regime, in
which the magnetic field magnitude was not allowed to vary much along
field lines. The proposed extension maintains the attractive features
of the original scheme (including the analytical Green's function,
which is critical for tractability) with minor modifications, while
allowing for completely general magnetic fields. The accuracy and
generality of the approach are demonstrated by numerical experiment
with an analytical manufactured solution.
\end{abstract}
\maketitle

\section{Introduction}

Recently, an asymptotic-preserving (AP) semi-Lagrangian scheme has
been proposed \citep{chacon2014asymptotic,chacon-jpc-24-imp_xport}
for the strongly anisotropic transport equation in magnetized plasmas
able to deal with arbitrary (including infinite \citep{del2011local,del2012parallel})
thermal-conductivity anisotropy ratios $\chi_{\parallel}/\chi_{\perp}$
(where $\chi_{\parallel}$ and $\chi_{\perp}$ are thermal conductivities
along and perpendicular to the magnetic field $\mathbf{B}$, respectively)
with numerical error independent of $\chi_{\parallel}/\chi_{\perp}$.\footnote{Theoretical estimates, experimental measurements, and modeling suggest
that the transport anisotropy in common tokamak reactors can reach
extremely high values $\chi_{\parallel}/\chi_{\perp}\sim10^{7}$ ---
$10^{10}$ \citep{braginskii1963transport,holzl2009determination,ren1998measuring,meskat2001analysis,snape2012influence,choi2014improved}.} The approach is based on a Green's function reformulation of the
parallel transport operator along magnetic field lines, which renders
the problem well conditioned with respect to the anisotropy ratio,
but at the expense of the formulation becoming integro-differential.
Operator-split \citep{chacon2014asymptotic} and fully implicit \citep{chacon-jpc-24-imp_xport}
solvers have been proposed to deal with the resulting formulation,
both unconditionally stable with respect to the timestep, and with
algorithmic performance independent of the anisotropy ratio. However,
these formulations were proposed for the so-called ``tokamak-ordering''
regime, in which the magnetic field compressibility {[}which we define
as $\nabla\cdot(\mathbf{B}/B)${]} is negligible.

In this study, we propose a practical generalization of the semi-Lagrangian
formulation for arbitrary magnetic-field compressibility that introduces
minimal changes to the tokamak-ordering one while retaining its attractive
features. Crucially, the new formulation is able to leverage the analytical
Green's functions that made the earlier approach tractable. Numerical
results with a manufactured solution \citep{salari2000code} demonstrate
the ability of the method to converge to the analytical solution for
arbitrary magnetic fields.

The rest of the paper is organized as follows. We review the tokamak-ordering
formulation in Sec. \ref{sec:formulation}. The generalization of
the method to arbitrary magnetic fields is introduced in Sec. \ref{sec:beyond-tokamak},
along with its favorable asymptotic properties. Details of the numerical
implementation are given in Sec. \ref{sec:num_imple}. Numerical results
demonstrating the properties of the scheme with a manufactured solution
are provided in Sec. \ref{sec:numerical_tests}, and we conclude in
Sec.~\ref{sec:conclusions}.

\section{Semi-Lagrangian algorithm in the tokamak-ordering regime}

\label{sec:formulation}We review briefly the large guide-field limit
(``tokamak-ordering'') semi-Lagrangian formulation. The anisotropic
transport equation, normalized to the perpendicular transport time
and length scales ($L_{\perp}$, $\tau_{\perp}=L_{\perp}^{2}/\chi_{\perp}$),
reads: 
\begin{equation}
\partial_{t}T-\frac{1}{\epsilon}\nabla_{\parallel}^{2}T=\nabla_{\perp}^{2}T+S\equiv S_{*},\label{diffusion_eq}
\end{equation}
where $T=T(t,\mathbf{x})$ is a temperature profile, $S=S(t,\mathbf{x})$
is a heat source, $S_{*}=\nabla_{\perp}^{2}T+S$ is a formal source
to the purely parallel transport equation, $\epsilon=\tau_{\parallel}/\tau_{\perp}=(L_{\parallel}^{2}/\chi_{\parallel})/(L_{\perp}^{2}/\chi_{\perp})$
is the ratio between parallel and perpendicular (to the magnetic field)
transport time scales, with $L_{\parallel},L_{\perp}$ being parallel
and perpendicular length scales, respectively. The thermal conductivity
along and perpendicular to the magnetic field ($\chi_{\parallel}$
and $\chi_{\perp}$) are assumed to be constants and uniform. The
differential operators along and perpendicular to the magnetic field
$\mathbf{B}=\mathbf{b}B$ are defined as:
\[
\nabla_{\parallel}^{2}=\nabla\cdot(\mathbf{b}\mathbf{b}\cdot\nabla),\quad\nabla_{\perp}^{2}=\nabla^{2}-\nabla_{\parallel}^{2}.
\]

The asymptotic limit equation (independent of $\epsilon$) of Eq.
\ref{diffusion_eq} is found by field-line averaging it as \citep{chacon2014asymptotic}:
\begin{equation}
\partial_{t}T_{\mathcal{N}}=\left\langle \nabla_{\perp}^{2}T+S\right\rangle =\left\langle \nabla^{2}T+S\right\rangle ,\label{eq:limit_problem_eq}
\end{equation}
where $T_{\mathcal{N}}=\left\langle T\right\rangle \in\mathcal{N}(\nabla_{\parallel})$
is in the null space of the parallel derivative operator, spanned
by constants along field lines. Here, $\left\langle A\right\rangle $
is the standard field-line-averaging operator, defined as:
\begin{equation}
\left\langle A\right\rangle =\frac{\int\frac{ds}{B(\hat{\mathbf{x}}(s;\mathbf{x}))}A(\hat{\mathbf{x}}(s;\mathbf{x}))}{\int ds/B(\hat{\mathbf{x}}(s;\mathbf{x}))}.\label{eq:field-line-avg}
\end{equation}
The field-line integration is performed along the magnetic field line
that passes through $\mathbf{x}$ and is parameterized by the arc
length~$s$,
\begin{align}
\frac{d\mathbf{\hat{x}}(s)}{ds}=\mathbf{b},\quad\mathbf{\hat{x}}(0)=\mathbf{x}.\label{eq:trace_mag}
\end{align}
The integral limits in Eq. \ref{eq:field-line-avg} can be finite
or infinite, depending on whether the field line is closed or open.
The asymptotic temperature field is given by:
\begin{equation}
T=T_{\mathcal{N}}+\mathcal{O}(\epsilon).\label{eq:asympt-T}
\end{equation}

In the tokamak-ordering regime (e.g., with a large toroidal field
$\mathbf{B}_{t}$ along which gradients are negligible), we assume
$\mathbf{B}=\mathbf{B}_{p}+\mathbf{B}_{t}$, $B_{t}\gg B_{p}$ , and
therefore: 
\[
\nabla\cdot\mathbf{b}=\mathbf{B}\cdot\nabla\frac{1}{B}\approx-\mathbf{B}_{p}\cdot\frac{\nabla B}{B_{t}^{2}}\ll1.
\]
Therefore, we can write:

\[
\nabla\cdot(\mathbf{b}\mathbf{b}\cdot\nabla)=\left(\cancelto{\approx0}{(\nabla\cdot\mathbf{b})}+(\mathbf{b}\cdot\nabla)\right)(\mathbf{b}\cdot\nabla)\approx(\mathbf{b}\cdot\nabla)^{2}=\frac{\partial^{2}}{\partial s^{2}}.
\]
The Green's function solution of the anisotropic transport equation
in the tokamak-ordering regime is given by \citep{chacon2014asymptotic}:
\begin{align}
T(t,\mathbf{x})=\mathcal{G}\left(T_{0};\mathbf{x},\frac{t}{\epsilon}\right)+\int_{0}^{t}dt'\mathcal{G}\left(S_{*};\mathbf{x},\frac{t-t'}{\epsilon}\right),\label{green:sol}
\end{align}
where 
\begin{equation}
\mathcal{G}\left(T_{0};\mathbf{x},t\right)=\int dsG(s,t)T_{0}\left(\mathbf{\hat{x}}(s,\mathbf{x})\right)\label{eq:prop}
\end{equation}
is the propagator of the homogeneous transport equation, $T(0,\mathbf{x})=T_{0}(\mathbf{x})=T_{0}$
is the initial condition, and $G(s,t)$ is the Green's function of
the diffusion equation, which in the case of infinite (perfectly confined)
magnetic field lines reads:
\begin{equation}
G(s,t)=\frac{1}{\sqrt{4\pi t}}\exp\left(-\frac{s^{2}}{4t}\right).\label{green:fun}
\end{equation}
The integration in Eq. \ref{eq:prop} is performed along the magnetic
field line, Eq. \ref{eq:trace_mag}. The semi-Lagrangian formulation
(Eq. \ref{green:sol}) features important properties, namely, that
$\mathcal{G}$ is the identity on $\mathcal{N}(\nabla_{\parallel})$,
and that $\lim_{t\to\infty}\mathcal{G}$ is the projector onto $\mathcal{N}(\nabla_{\parallel})$.
These play a central role in controlling numerical pollution, and
ensuring the asymptotic preserving properties of any numerical method
constructed based on Eq. \ref{green:sol} when $\epsilon\rightarrow0$
\citep{chacon2014asymptotic}.

Equation \ref{green:sol} can be effectively discretized in time by
approximating $S_{*}$ as a constant in time \citep{chacon2014asymptotic,chacon-jpc-24-imp_xport}.
For first-order implicit Backward Euler (or Backward Differentiation
Formula of order unity, BDF1), Eq. \ref{green:sol} gives:
\begin{align}
T(\mathbf{x})^{n+1}=\mathcal{G}\left(T^{n};\mathbf{x},\frac{\Delta t}{\epsilon}\right)+\Delta t\mathcal{P}\left(S_{*}^{n+1};\mathbf{x},\frac{\Delta t}{\epsilon}\right)+\mathcal{O}(\Delta t^{2}),\label{bdf1}
\end{align}
where 
\begin{align}
\mathcal{P}\left(S_{*}^{n+1};\mathbf{x},\frac{\Delta t}{\epsilon}\right)=\int_{-\infty}^{+\infty}ds\mathcal{U}\left(s,\frac{\Delta t}{\epsilon}\right)S_{*}^{n+1}\left(\mathbf{\hat{x}}(s,\mathbf{x})\right),\label{prop_def}
\end{align}
and 
\begin{align}
\mathcal{U}\left(s,\tau\right) & =\frac{1}{\sqrt{\tau}}\left(\frac{e^{-s^{2}/4\tau}}{\sqrt{\pi}}-\frac{|s|}{2\sqrt{\tau}}\mathrm{erfc}\left(\frac{|s|}{2\sqrt{\tau}}\right)\right).
\end{align}
Equation \ref{bdf1} can be used to construct higher-order BDF formulas
\citep{chacon2014asymptotic,chacon-jpc-24-imp_xport}. It automatically
yields the correct asymptotic limit when taking the $\epsilon\rightarrow0$
limit by using the projection property of the propagators $\mathcal{G}$
and $\mathcal{P}$ \citep{chacon2014asymptotic}, to find the time-discrete
limit equation:
\begin{equation}
T_{\mathcal{N}}^{n+1}=T_{\mathcal{N}}^{n}+\Delta t\left\langle \nabla_{\perp}^{2}T+S\right\rangle _{s}=T_{\mathcal{N}}^{n}+\Delta t\left\langle \nabla^{2}T+S\right\rangle _{s},\label{eq:limit_problem_discrete}
\end{equation}
where $\left\langle A\right\rangle _{s}$ is the tokamak-ordering
field-line-averaging operator:
\begin{equation}
\left\langle A\right\rangle _{s}=\frac{\int dsA(\hat{\mathbf{x}}(s;\mathbf{x}))}{\int ds}.\label{eq:to-field-line-avg}
\end{equation}

\section{Generalization to arbitrary magnetic fields}

\label{sec:beyond-tokamak}We generalize next the semi-Lagrangian
transport approach to $\nabla\cdot\mathbf{b}\neq0$. For arbitrary
magnetic fields, $\nabla\cdot(\mathbf{b}\mathbf{b}\cdot\nabla T)=B\partial_{s}(1/B\,\partial_{s}T)$,
and hence Eq. \ref{diffusion_eq} reads:
\begin{equation}
\partial_{t}T-\frac{B}{\epsilon}\partial_{s}\left(\frac{1}{B}\partial_{s}T\right)=\nabla_{\perp}^{2}T+S,\label{eq:parallel_xport_general_B}
\end{equation}
with $B=B[\mathbf{x}(s)]$ the magnitude of the magnetic field along
the field line. At magnetic nulls (where $B=0$), Eq. \ref{eq:parallel_xport_general_B}
becomes isotropic:
\begin{equation}
\partial_{t}T=\nabla^{2}T+S.\label{eq:xport_B=00003D0}
\end{equation}

An analytical expression for the Green's function of the parallel
transport term in Eq. \ref{eq:parallel_xport_general_B} does not
exist for arbitrary $B(s)$ profiles. However, analytical progress
is possible by rewriting Eq. \ref{eq:parallel_xport_general_B} in
terms of a new variable $\lambda$, related to the arc-length by:
\begin{equation}
\frac{d\lambda}{ds}=B(s).\label{eq:lambda-def}
\end{equation}
There results:
\begin{equation}
\frac{1}{B^{2}}\partial_{t}T-\frac{1}{\epsilon}\partial_{\lambda}^{2}T=\frac{\nabla_{\perp}^{2}T+S}{B^{2}}.\label{eq:xport-eq-lambda}
\end{equation}
This equation is in principle ill-posed for $B$ strictly equal to
zero, but this is never the case along a given magnetic field line.
It can only occur at magnetic nulls, and there is a simple prescription
within our framework to deal with this special case, which we outline
later in this section.

Adding and subtracting $\beta\partial_{t}T$, with $\beta(\mathbf{x})>0$
a \textcolor{black}{field to be defined precisely later, but considered
for the time being a constant}\textcolor{black}{\emph{ }}\textcolor{black}{along
the field line passing through $\mathbf{x}$}\textcolor{black}{\emph{
}}\textcolor{black}{over the kernel integration domain}, there results:
\begin{equation}
\beta\partial_{t}T-\frac{1}{\epsilon}\partial_{\lambda}^{2}T=\frac{\nabla_{\perp}^{2}T+S-\partial_{t}T}{B^{2}}+\beta\partial_{t}T.\label{eq:reform_t_eq}
\end{equation}
This equation admits a formal Lagrangian treatment by considering
a new formal source:
\begin{equation}
S^{*}=\frac{\nabla_{\perp}^{2}T+S-\partial_{t}T}{B^{2}}+\beta\partial_{t}T,\label{eq:formal_source_beta}
\end{equation}
which leads to the following first-order BDF1 scheme (assuming $\Delta t\,\partial_{t}T\ll T$
and $\Delta t\,\partial_{t}S\ll S$ \textcolor{black}{for accuracy purposes}):
\begin{equation}
T^{n+1}(\mathbf{x})=\mathcal{G}_{{\color{black}\lambda}}\left(T^{n};\mathbf{x},\frac{\Delta t}{\beta\epsilon}\right)+\frac{\Delta t}{\beta}\,\mathcal{P}_{\lambda}\left(\left.\frac{\nabla_{\perp}^{2}T+S}{B^{2}}\right|^{n+1}+\left(\beta-\frac{1}{B^{2}}\right){\color{black}\frac{T^{n+1}-T^{n}}{\Delta t}};\mathbf{x},\frac{\Delta t}{\beta\epsilon}\right)+\mathcal{O}(\Delta t^{2}).\label{eq:formal_solution_general_B}
\end{equation}
\textcolor{black}{We comment on the numerical implications of the presence
of the temporal derivative of the temperature field in the argument
of the kernel integral $\mathcal{P}_{\lambda}$ in Eq. \ref{eq:formal_solution_general_B}
due to the formal source in Eq. \ref{eq:formal_source_beta} later
in this study.}\textcolor{black}{{} Note that a second-order BDF scheme
can be formulated as well, following the prescription provided in
Ref. \citep{chacon2014asymptotic}.} The subscript $\lambda$ indicates
that the field-line integrals are now performed in the variable $\lambda$
instead of the arc-length $s$, i.e.:
\begin{eqnarray*}
\mathcal{G}_{\lambda}(A;\mathbf{x},\tau) & = & \int_{-\infty}^{+\infty}d\lambda'\,A[\mathbf{x}_{\lambda}(\lambda';\mathbf{x})]G(\lambda-\lambda',\tau),\\
\mathcal{P}_{\lambda}(A;\mathbf{x},\tau) & = & \int_{-\infty}^{\infty}d\lambda'\,A[\mathbf{x}_{\lambda}(\lambda';\mathbf{x})]\,\mathcal{U}(\lambda-\lambda',\tau).
\end{eqnarray*}
Here, $\mathbf{x}_{\lambda}(\lambda';\mathbf{x})$ solves the modified
magnetic field ODE:
\begin{equation}
\frac{d\mathbf{x}_{\lambda}}{d\lambda'}=\frac{\mathbf{b}}{B}\,\,,\,\,\mathbf{x}_{\lambda}(\lambda'=\lambda;\mathbf{x})=\mathbf{x}.\label{eq:B-field-ODE-lambda}
\end{equation}
As with the arc-length integrals, one can define a field-line average
in terms of the variable $\lambda,$ which annihilates the $\partial_{\lambda}^{2}$
in Eq. \ref{eq:xport-eq-lambda}, as:
\begin{equation}
\left\langle A\right\rangle _{\lambda}=\frac{\int d\lambda A(\mathbf{x}_{\lambda}(\lambda;\mathbf{x}))}{\int d\lambda}=\frac{\int A(\hat{\mathbf{x}}(s;\mathbf{x}))\,B(\hat{\mathbf{x}}(s;\mathbf{x}))ds}{\int B(\hat{\mathbf{x}}(s;\mathbf{x}))ds},\label{eq:lambda-average}
\end{equation}
where we have used Eq. \ref{eq:lambda-def}. The propagators $\mathcal{G}_{\lambda}$
and $\mathcal{P}_{\lambda}$ and the $\lambda$-averaging operator
$\left\langle \cdots\right\rangle _{\lambda}$ satisfy several important
properties (see Ref. \citep{chacon2014asymptotic} for proofs):
\begin{enumerate}[label=\roman{enumi}., ref=(\roman{enumi})]
\item As stated above, $\left\langle \partial_{\lambda}^{2}A\right\rangle _{\lambda}=0$,
i.e., the $\lambda$-average annihilates the parallel transport operator
in Eq. \ref{eq:xport-eq-lambda}.
\item \label{enu:lambda-aver-pprty-2}$\left\langle A_{\mathcal{N}}\right\rangle _{\lambda}=A_{\mathcal{N}}$,
i.e., the $\lambda$-averaging operator is the identity on $\mathcal{N}(\nabla_{\parallel})$.
This property simply states that the average of a constant is the
same constant.
\item \label{enu:P-average}$\left\langle \mathcal{P}_{\lambda}(A;\mathbf{x},\tau)\right\rangle _{\lambda}=\left\langle A\right\rangle _{\lambda}$
and $\left\langle \mathcal{G}(A;\mathbf{x},\tau)\right\rangle _{\lambda}=\left\langle A\right\rangle _{\lambda}$.
\item \label{enu:P-project}$\mathcal{G}_{\lambda}(A;\mathbf{x},\tau\rightarrow\infty)\rightarrow\left\langle A\right\rangle _{\lambda}=A_{\mathcal{N}}\in\mathcal{N}(\nabla_{\parallel})$,
$\mathcal{P}_{\lambda}(A;\mathbf{x},\tau\rightarrow\infty)\rightarrow\left\langle A\right\rangle _{\lambda}=A_{\mathcal{N}}\in\mathcal{N}(\nabla_{\parallel})$,
i.e., the propagators become projectors to the null space when $\tau\rightarrow\infty$.
\item \label{enu:P-asymp}The $\mathcal{P}_{\lambda}$-propagator satisfies
the following asymptotic properties:
\begin{equation}
\mathcal{P}_{\lambda}(A;\mathbf{x},\tau)=\begin{cases}
\left\langle A\right\rangle _{\lambda}+\mathcal{O}(\tau^{-1}), & \tau\gg1\\
A(\mathbf{x})[1+\mathcal{O}(\tau)], & \tau\ll1
\end{cases}.\label{eq:P-asymp-properties}
\end{equation}
\item \label{enu:lambda-aver-pprty-1}$\left\langle \frac{A}{B^{2}}\right\rangle _{\lambda}=\frac{\int A/Bds}{\int Bds}=\left\langle A\right\rangle \frac{\int ds/B}{\int Bds}=\frac{\left\langle A\right\rangle }{\left\langle B^{2}\right\rangle }$,
with $\left\langle A\right\rangle $ defined in Eq. \ref{eq:field-line-avg}.
This property is crucial to ensure that Eq. \ref{eq:xport-eq-lambda}
features the same limit equation as Eq.~\ref{eq:parallel_xport_general_B}
when $\epsilon\rightarrow0$.
\item \label{enu:lambda-AP-limit}$\left\langle T\right\rangle _{\lambda}\rightarrow\left\langle T\right\rangle =T_{\mathcal{N}}$
as $\epsilon\rightarrow0$, i.e., both averaging operators give the
same null space component of the solution. This follows from $\left\langle T\right\rangle _{\lambda}=\left\langle \left\langle T\right\rangle +\mathcal{O}\left(\epsilon\right)\right\rangle _{\lambda}=T_{\mathcal{N}}+\mathcal{O}\left(\epsilon\right)$.
\end{enumerate}
It can be readily shown that the correct limit equation (Eq. \ref{eq:limit_problem_discrete})
follows in the $\epsilon\rightarrow0$ limit when taking the $\lambda$-average
of Eq. \ref{eq:formal_solution_general_B}. Using properties \ref{enu:P-average}
and \ref{enu:lambda-aver-pprty-1} of the $\lambda$-average,\textcolor{black}{{}
and that $\beta$ becomes a field-line constant when $\epsilon\rightarrow0$
as shown below}, we find from Eq. \ref{eq:formal_solution_general_B}
that:
\[
\left\langle T^{n+1}\right\rangle _{\lambda}=\left\langle T^{n}\right\rangle _{\lambda}+\frac{\Delta t}{\beta}\left[\left.\frac{\left\langle \nabla_{\perp}^{2}T+S\right\rangle }{\left\langle B^{2}\right\rangle }\right|^{n+1}+\left(\beta-\frac{1}{\left\langle B^{2}\right\rangle }\right)\left\langle \frac{T^{n+1}-T^{n}}{\Delta t}\right\rangle _{\lambda}\right],
\]
which immediately leads to Eq. \ref{eq:limit_problem_discrete}.

We motivate next the particular choice of $\beta$ used in this work
to produce a robust, AP, and convergent numerical scheme.

\subsection{Definition of  $\beta(\mathbf{x})$}

\label{subsec:beta-def}

To guide the choice of $\beta(\mathbf{x})$, we reformulate Eq. \ref{eq:formal_solution_general_B}
in terms of $\Delta T=T^{n+1}-T^{n}$ as:
\begin{equation}
\Delta T-\mathcal{P}_{\lambda}\left(\left[1-\frac{1}{\beta B^{2}}\right]\Delta T;\mathbf{x},\frac{\Delta t}{\beta\epsilon}\right)=\mathcal{G}_{\lambda}\left(T^{n};\mathbf{x},\frac{\Delta t}{\beta\epsilon}\right)-T^{n}+\frac{\Delta t}{\beta}\,\mathcal{P}_{\lambda}\left(\left.\frac{\nabla_{\perp}^{2}T+S}{B^{2}}\right|^{n+1};\mathbf{x},\frac{\Delta t}{\beta\epsilon}\right).\label{eq:general-B-delta-form}
\end{equation}
\textcolor{black}{The right-hand side of this equation is similar to
the tokamak-ordering one (except for the presence of $\beta$ and
$1/B^{2}$ terms), and could in principle be computed with an operator-split
approach as proposed in Ref. \citep{chacon2014asymptotic}. However,
the left-hand side of Eq. \ref{eq:general-B-delta-form} contains
an integral operator on $\Delta T$, which cannot be operator-split
in the same fashion because, as our analysis in \ref{app:os-error}
shows, the resulting scheme is not convergent with either $\Delta t$
or $\epsilon$. This is nor surprising, since a strict numerical balance
must be struck between the temporal derivative terms in the reformulated
temperature equation (Eq. \ref{eq:reform_t_eq}) for the formulation
to be equivalent to the original one. It therefore needs to be iterated
for accuracy. Consequently, we do not consider the operator-split
formulation further in this study, and we focus on the fully implicit
discretization of Eq. \ref{eq:formal_solution_general_B} (i.e., Eq.
\ref{eq:general-B-delta-form}).}

\textcolor{black}{The conditioning of the integral operator of the
left-hand side of Eq. \ref{eq:general-B-delta-form} can be improved
by choosing $\beta$ appropriately.} In this study, we choose $\beta$
such that the contribution of the $\mathcal{P}_{\lambda}$ integral
in the left hand side is zero for the null space component of $\Delta T$,
$\left\langle \Delta T\right\rangle $, i.e., we require:
\begin{equation}
\mathcal{P}_{\lambda}\left(1-\frac{1}{\beta B^{2}};\mathbf{x},\frac{\Delta t}{\beta\epsilon}\right)=0\Rightarrow\beta(\mathbf{x})=\mathcal{P}_{\lambda}\left(\frac{1}{B^{2}};\mathbf{x},\frac{\Delta t}{\beta\epsilon}\right),\label{eq:beta-def-1}
\end{equation}
\textcolor{black}{where in the last step we have used that $\beta$
is by assumption a constant along the field line on the integration
domain of the kernel integral. }This choice restricts the contribution
of $\mathcal{P}_{\lambda}$ in the left hand side of Eq. \ref{eq:general-B-delta-form}
to be of $\mathcal{O}(\epsilon)$ (very small except at boundary layers
such as island separatrices) and therefore generally small vs. $\Delta T$,
facilitating a faster convergence of Eq. \ref{eq:general-B-delta-form}
in an iterative context. Equation \ref{eq:beta-def-1} is a nonlinear
definition of $\beta$ that is only a function of the magnetic field
topology (and not of the temperature field), and also requires iteration.
We discuss a practical way of finding $\beta$ in Sec.~\ref{sec:num_imple}.

\textcolor{black}{We show next that the choice for $\beta$ in Eq.
\ref{eq:beta-def-1}}\textcolor{black}{{} }enforces the correct asymptotic
limits. From Eq. \ref{eq:beta-def-1}, and the asymptotic properties
of the $\mathcal{P}$ propagator (Eq. \ref{eq:P-asymp-properties}
in property \ref{enu:P-asymp} of the $\mathcal{P}_{\lambda}$-propagator),
$\beta$ at a given spatial point $\mathbf{x}$ has the following
asymptotic properties:
\[
\beta(\mathbf{x})=\begin{cases}
\left\langle \frac{1}{B^{2}}\right\rangle _{\lambda}+\mathcal{O}\left(\frac{\beta\epsilon}{\Delta t}\right), & \Delta t\gg\beta\epsilon\\
\frac{1}{B^{2}(\mathbf{x})}\left[1+\mathcal{O}\left(\frac{\Delta t}{\beta\epsilon}\right)\right], & \Delta t\ll\beta\epsilon
\end{cases}.
\]
By property \ref{enu:lambda-aver-pprty-1} of the $\lambda$-averaging
operator, we have:
\[
\left\langle \frac{1}{B^{2}}\right\rangle _{\lambda}=\frac{1}{\left\langle B^{2}\right\rangle },
\]
and therefore:
\begin{equation}
\beta(\mathbf{x})=\begin{cases}
\frac{1}{\left\langle B^{2}\right\rangle }\left[1+\mathcal{O}\left(\frac{\epsilon}{\Delta t}\right)\right], & \left\langle B^{2}\right\rangle \Delta t\gg\epsilon\\
\frac{1}{B^{2}(\mathbf{x})}\left[1+\mathcal{O}\left(\frac{B^{2}\Delta t}{\epsilon}\right)\right], & B^{2}\Delta t\ll\epsilon
\end{cases}.\label{eq:beta-asymp-properties}
\end{equation}
It is noteworthy that $\beta(\mathbf{x})\rightarrow1/B^{2}(\mathbf{x})$
(or, alternatively, $\beta B^{2}\rightarrow1$) when $B\rightarrow0$
for arbitrary $\epsilon>0$. This ensures that Eq. \ref{eq:general-B-delta-form}
is well posed when $B$ becomes arbitrarily small, as we shall see.
An additional advantage of the definition of $\beta$ in Eq. \ref{eq:beta-def-1}
is that it automatically gives $\beta=1/B_{0}^{2}$ when $B=B_{0}$
is a constant along magnetic fields. This, in turn, ensures that the
time-discrete Eq. \ref{eq:general-B-delta-form} reverts back to Eq.
\ref{bdf1}.

In the next section, we investigate the asymptotic properties of Eq.
\ref{eq:general-B-delta-form}.

\subsection{Asymptotic properties of Eq. \ref{eq:formal_solution_general_B}}

We are interested in studying three distinct limits: $\Delta t\rightarrow0$
(consistency), $\epsilon\rightarrow0$ (asymptotic preservation),
and $B\rightarrow0$ (regularity at magnetic nulls). Consistency can
be readily demonstrated by taking the $\Delta t\rightarrow0$ limit
in Eq. \ref{eq:general-B-delta-form}, and using the asymptotic properties
of $\mathcal{P}_{\lambda}$, $\mathcal{G}_{\lambda}$, and of $\beta\rightarrow1/B^{2}$,
to find:
\[
\Delta T=\frac{B^{2}\Delta t}{\epsilon}\partial_{\lambda}^{2}T^{n}+\Delta t\left[\nabla_{\perp}^{2}T+S\right]^{n+1}+\mathcal{O}(\Delta t^{2})\stackrel{\Delta t\rightarrow0}{\longrightarrow}\partial_{t}T-\frac{B}{\epsilon}\partial_{s}\left(\frac{1}{B}\partial_{s}T\right)=\nabla_{\perp}^{2}T+S.
\]
Similarly, regularity at magnetic nulls can be demonstrated from Eq.
\ref{eq:general-B-delta-form} by using that, for arbitrary $\Delta t$
and $\epsilon$, as $B\rightarrow0$, $\Delta t/\beta\epsilon\sim B^{2}\rightarrow0$
and $\beta B^{2}\rightarrow1$, and therefore:
\[
\Delta T=\Delta t\left[\nabla_{\perp}^{2}T+S\right]^{n+1}+\mathcal{O}(B^{2}\Delta t/\epsilon)\stackrel{B\rightarrow0}{\longrightarrow}\Delta t\left[\nabla^{2}T+S\right]^{n+1}.
\]

We have already shown in Sec. \ref{sec:beyond-tokamak} that Eq. \ref{eq:general-B-delta-form}
leads to the correct limit problem when $\lambda$-averaging it. Here,
we show that it is AP by taking its $\epsilon\rightarrow0$ limit.
Using the asymptotic property \ref{enu:P-project} of $\mathcal{G}_{\lambda}$,
$\mathcal{P}_{\lambda}$, we find:
\[
\Delta T-\left\langle \left[1-\frac{1}{\beta B^{2}}\right]\Delta T\right\rangle _{\lambda}-\left\langle T^{n}\right\rangle _{\lambda}+T^{n}-\Delta t\left\langle \frac{\left[\nabla_{\perp}^{2}T+S\right]^{n+1}}{\beta B^{2}}\right\rangle _{\lambda}=0.
\]
Introducing $\Delta T=T^{n+1}-T^{n}$, and using properties \ref{enu:lambda-aver-pprty-1}
and \ref{enu:lambda-AP-limit} of the $\lambda$-average and the asymptotic
properties of $\beta$ (Eq. \ref{eq:beta-asymp-properties}), we find:
\[
T^{n+1}=\left\langle T^{n}\right\rangle +\Delta t\left\langle \left[\nabla^{2}T+S\right]^{n+1}\right\rangle ,
\]
which is the time-discrete formulation of the limit equation. Therefore,
$T^{n+1}\rightarrow\left\langle T^{n+1}\right\rangle $ as $\epsilon\rightarrow0$,
and the formulation in Eq. \ref{eq:general-B-delta-form} is AP.

Therefore, we conclude that Eq. \ref{eq:general-B-delta-form} (and
therefore Eq. \ref{eq:formal_solution_general_B}) satisfies the following
limits:
\[
\mathrm{Eq.\,(\ref{eq:general-B-delta-form})}\rightarrow\begin{cases}
\partial_{t}T-\frac{B}{\epsilon}\partial_{s}\left(\frac{1}{B}\partial_{s}T\right)=\nabla_{\perp}^{2}T+S, & \Delta t\rightarrow0,\,\forall\epsilon,B\neq0\\
\left\langle T^{n+1}\right\rangle =\left\langle T^{n}\right\rangle +\Delta t\left\langle \left[\nabla^{2}T+S\right]^{n+1}\right\rangle , & \epsilon\rightarrow0,\,\forall\Delta t,B\neq0\\
T^{n+1}=T^{n}+\Delta t\left[\nabla^{2}T+S\right]^{n+1}, & B\rightarrow0,\,\forall\Delta t,\epsilon
\end{cases}.
\]
The first limit establishes consistency, the second its AP nature,
and the third ensures regularity of the formulation at magnetic nulls.

\section{Numerical implementation details}

\label{sec:num_imple}

\textcolor{black}{The discrete fields $T^{n+1}$, $T^{n}$, $S^{n+1}$,
$\beta$ and the magnetic field are collocated on a computational
grid. We employ the same fourth-order discretization for all spatial
operators as used in Refs. \citep{chacon2014asymptotic,chacon-jpc-24-imp_xport}.
The Lagrangian integrals in the operators $\mathcal{G}_{\lambda}$
and $\mathcal{P}_{\lambda}$ require the reconstruction (by interpolation)
of these discrete fields over the whole domain to evaluate them at
arbitrary points along magnetic field orbits. This is done in this
study with global, arbitrary-order splines, but we have also implemented
and tested second-order B-spline-based positivity-preserving local-stencil
interpolations \citep{chacon2022asymptotic}. Additional numerical
implementation details for the field-line integrals are provided in
Refs. \citep{chacon2014asymptotic,chacon-jpc-24-imp_xport}. }

We employ the GMRES-based implicit algorithm of Ref. \citep{chacon-jpc-24-imp_xport}
as the basis for all the simulations presented here. This straightforwardly
deals with the iteration needed to solve Eq. \ref{eq:general-B-delta-form},
which is well conditioned by construction by the choice of $\beta$.

For the tests presented here, the magnetic field is static in time,
and therefore we only need to solve for $\beta$ once at the beginning
of the simulation. The field $\beta$ must be determined\textcolor{black}{{}
at every mesh cell }from the nonlinear equation Eq. \ref{eq:beta-def-1}.
We use a Picard iteration for this. We begin the iteration by providing
an initial value of $\beta$, $\beta^{0}$, at every mesh point as:
\[
\beta^{0}=\frac{1}{B^{2}},
\]
where $B$ is the cell magnetic field magnitude. This choice ensures
that $\beta B^{2}\rightarrow1$ when $B\rightarrow0$, which is important
to control numerical errors in this limit. In subsequent iterations,
we update $\beta$ according to a simple Picard iterative prescription:
\begin{equation}
\beta^{k+1}(\mathbf{x})=\mathcal{P}_{\lambda}\left(\frac{1}{B^{2}};\mathbf{x},\frac{\Delta t}{\beta^{k}\epsilon}\right).\label{eq:beta-iteration}
\end{equation}
When $\Delta t\gg\beta\epsilon$, this iteration is expected to be
contractive (and therefore convergent) due to the strongly smoothing
property of the integral propagator. In the opposite limit, $\Delta t\ll\beta\epsilon$,
we have (Eq. \ref{eq:P-asymp-properties}):
\[
\mathcal{P}_{\lambda}\left(\frac{1}{B^{2}};\mathbf{x},\frac{\Delta t}{\beta\epsilon}\right)=\frac{1}{B(\mathbf{x})^{2}}\left[1+\mathcal{O}\left(\frac{\Delta t}{\beta\epsilon}\right)\right],
\]
and therefore our initial guess is almost exact. We converge this
iteration to a relative tolerance of $10^{-5}$. In practice, we find
$\beta$ converges in at most four or five iterations. In a time-varying
magnetic-field context, we plan to use a single iteration of Eq. \ref{eq:beta-iteration}
to update $\beta$ in time, i.e., we will prescribe $k=n$ and $n+1=k+1$.

\section{Numerical tests}

\label{sec:numerical_tests}

We consider a manufactured solution in two dimensions for our numerical
tests with homogeneous Dirichlet boundary conditions in $x$ and periodic
boundary conditions in $y$, wit\textcolor{black}{h $x,y\in[0,1]^{2}$
and $L_{\parallel}=L_{\perp}$. The} magnetic field is of the form:
\begin{equation}
{\bf B}=\mathbf{z}\times\nabla\psi+B_{0}\mathbf{z},\label{magn_field}
\end{equation}
with $B_{0}$ a constant. The magnetic topology is determined by $\psi$
and the \textbf{b}-compressibility by $B_{0}$, since for this choice
$\nabla\cdot\mathbf{b}=\mathcal{O}\left(\frac{B_{p}}{L_{\perp}B_{0}}\right)^{3}$.
Increasing $B_{0}$ will result in convergence to the tokamak-ordering
solution. We consider a flux function $\psi$ of the form:
\begin{equation}
\psi=x+\delta\sin(2\pi x)\cos(2\pi y).
\end{equation}
The corresponding magnetic field components are $B_{x}=-\frac{\partial\psi}{\partial y}=2\pi\delta\sin(2\pi x)\sin(2\pi y)$,
$B_{y}=\frac{\partial\psi}{\partial x}=1+2\pi\delta\cos(2\pi x)\cos(2\pi y)$
and $B_{z}=B_{0}$. Note that $B_{p}=\sqrt{B_{x}^{2}+B_{y}^{2}}\sim\mathcal{O}(1)$
on average. 
\begin{figure*}
\centering{}\includegraphics[width=0.48\columnwidth]{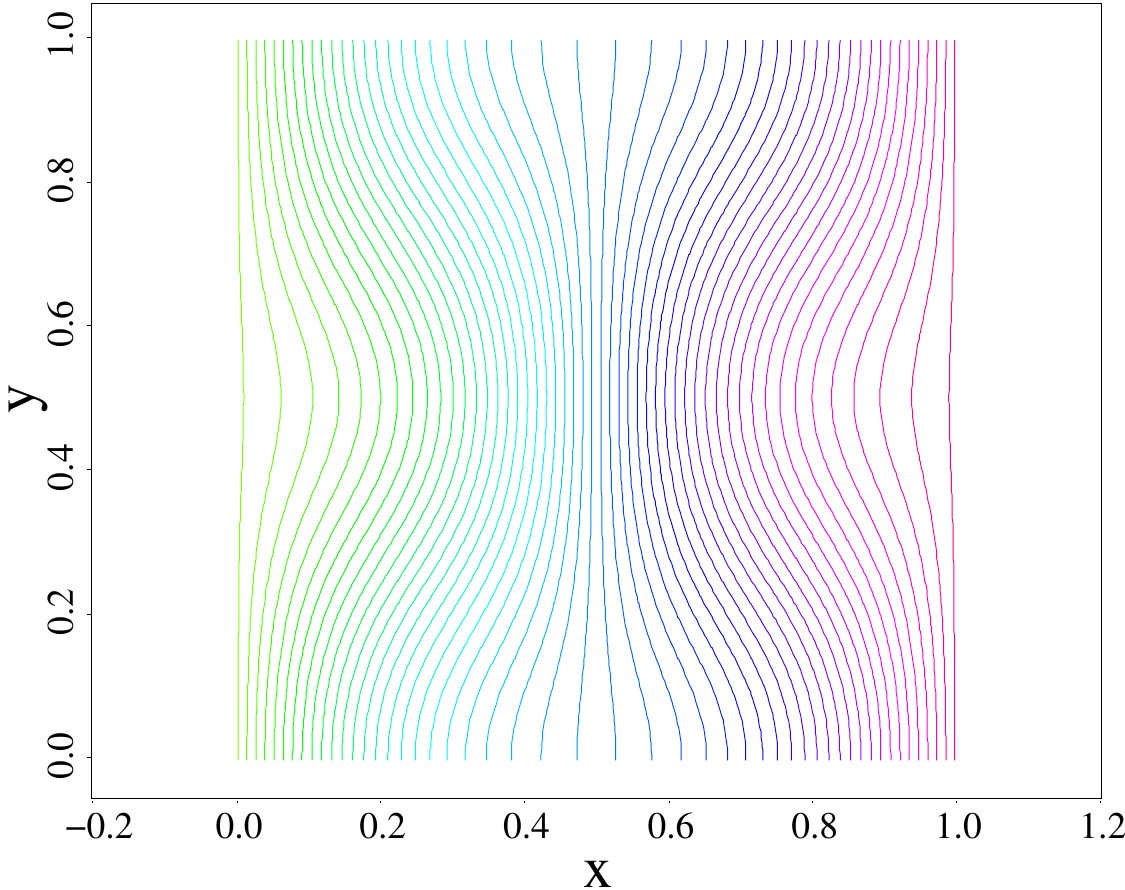}
\includegraphics[width=0.48\columnwidth]{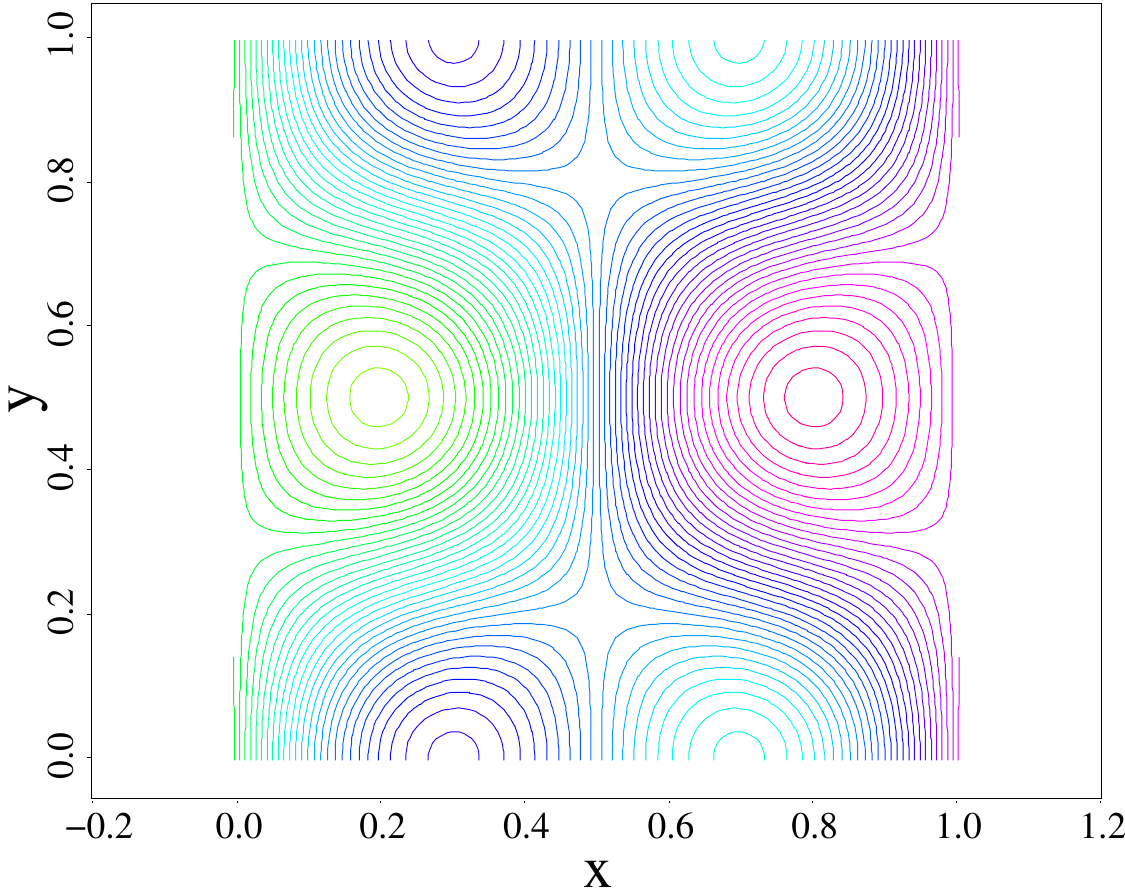}
\caption{Left panel: flux function with $\delta=0.1$. Right panel: flux function
with $\delta=0.5$.\label{fig:flux_function}}
\end{figure*}
As shown in Fig.~\ref{fig:flux_function}, $\delta$ controls the
magnetic field topology: for $\delta=0.1$ the magnetic field is simply
connected, whereas for $\delta=0.5$ it features magnetic islands
with separatrices, X-points, and O-points. The latter are special
points satisfying $\nabla\psi=\mathbf{0}$, and are magnetic nulls
when $B_{0}=0$.

We consider a steady-state temperature field of the form: 
\begin{equation}
T_{\infty}(x,y)=T_{\mathcal{N}}(\psi)+\epsilon\:\tilde{T}(x,y).\label{eq:T_ansatz}
\end{equation}
Here, $\epsilon$ is the dimensionless anisotropy ratio, as defined
above. This ansatz is consistent with the asymptotic properties of
the temperature field as $\epsilon\rightarrow0$, Eq. \ref{eq:asympt-T}.
In particular, we choose $T_{\mathcal{N}}(\psi)=\psi$ and $\tilde{T}=\cos(4\pi y)\sin(4\psi x)$.
This choice is consistent with the boundary conditions, and guarantees
that $\nabla_{\parallel}^{2}\tilde{T}\not=0$.

The manufactured solution source is found by inserting the prescribed
steady-state temperature (Eq. \ref{eq:T_ansatz}) into the steady-state
transport equation, to find:
\begin{equation}
S=-(1-\epsilon)\nabla_{\parallel}^{2}\tilde{T}-\epsilon\nabla^{2}\tilde{T}-\nabla^{2}T_{\mathcal{N}}(\psi).
\end{equation}
The explicit evaluation of the source is given in \ref{sec:Manufactured-source}.
We note that this particular source is singular when magnetic nulls
with $B=0$ are present (e.g., for $B_{0}=0$ and $\delta=0.5$).
This will have implications for spatial convergence, as we shall see.

We have performed tests for $\epsilon=10^{-2},10^{-4},10^{-6}$, $\delta=0.1,0.5$,
and various choices for $B_{0}$. Simulations are initialized using
the analytical steady-state solution, and are run to numerical steady
state to assess spatial accuracy (note that the temporal convergence
properties of the semi-Lagrangian scheme have been assessed elsewhere
\citep{chacon2014asymptotic,chacon-jpc-24-imp_xport}). For all simulations,
we use a BDF1 timestep $\Delta t=10^{-3}$ \textcolor{black}{(normalized
to $\chi_{\perp}$, which is comparable or larger to what would be
typical in an evolving B-field context, e.g., magnetohydrodynamics)}
and run them until the error saturates, usually for 100 to 200 steps.
Absolute errors are computed as an $\ell_{2}$-norm of the difference
between the numerical solution and the analytical one, Eq. \ref{eq:T_ansatz},
\begin{equation}
\mathrm{Error}=\sum_{i,j}\Delta x_{i}\Delta y_{j}(T_{ij}-T_{\infty}(x_{i},y_{j}))^{2}.\label{eq:l2-norm}
\end{equation}

\subsection{Convergence study of the tokamak-ordering formulation}

\begin{figure*}
\centering{}\includegraphics[width=0.48\columnwidth]{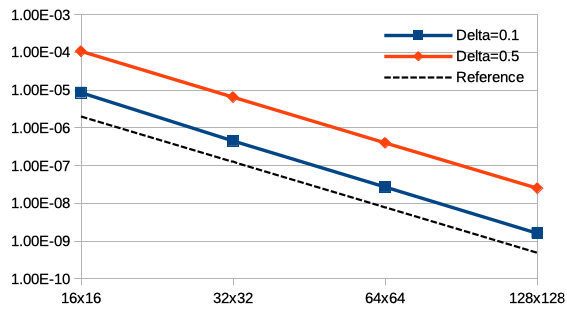}
\caption{Mesh convergence study for the tokamak-ordering formulation with a
null-space solution ($\tilde{T}=0$) for $\epsilon=10^{-2}$ with
$B_{0}=0$ and $\delta=0.1,0.5$. Fourth-order spatial convergence
is achieved, as expected. \label{fig:null-space-sol}}
\end{figure*}

We begin with a mesh convergence study of the tokamak-ordering formulation
with the null-space solution, i.e., $\tilde{T}=0$, which is insensitive
to the magnetic-field topology (since the $\nabla_{\parallel}$ operator
is annihilated). This test is intended to confirm the expected spatial
discretization accuracy (fourth-order), and has been performed with
$\delta=0.1,0.5$ and with $B_{0}=0$ (i.e., no guide field). The
results are shown in Fig \ref{fig:null-space-sol}, and confirm the
expected fourth-order accuracy scaling. 

\begin{figure*}
\centering{}\includegraphics[width=0.48\columnwidth]{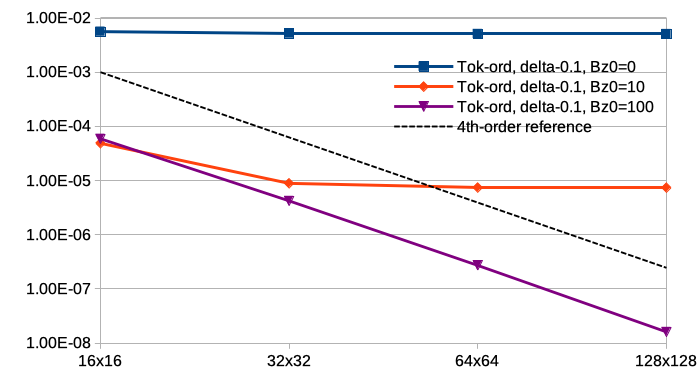}
\includegraphics[width=0.48\columnwidth]{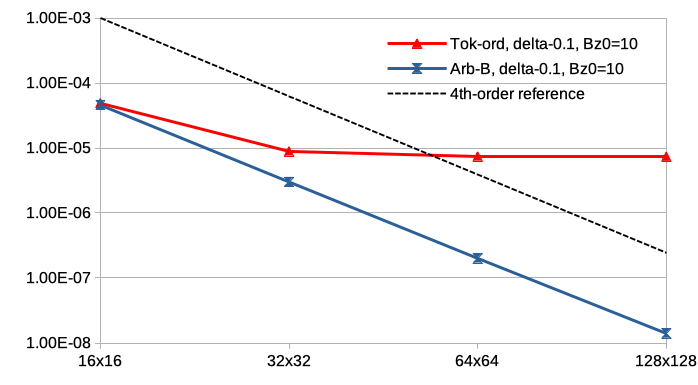}\caption{Mesh convergence study for the full solution with $\epsilon=10^{-2}$
and $\delta=0.1$. Left: tokamak-ordering formulation for three different
guide fields, $B_{0}=0,10,100$, demonstrating convergence of the
tokamak-ordering solution to the analytical one as the guide-field
$B_{0}$ increases. Right: Comparison of errors between tokamak-ordering
and arbitrary-B formulation for $B_{0}=10$, demonstrating the ability
of the arbitrary-B formulation to converge where the tokamak-ordering
one fails. \label{fig:tokamak_ordering_vs_arbB}}
\end{figure*}

Next, we confirm that the tokamak-ordering formulation converges
to the true analytical solution as the guide field $B_{0}$ is increased.
We consider the case of $\delta=0.1$ and $\epsilon=10^{-2}$. The
results are shown in Fig. \ref{fig:tokamak_ordering_vs_arbB}-left.
 As expected, with $B_{0}=0$ there is a strong \textbf{b}-compressibility
and the tokamak ordering approximation is unable to reproduce the
temperature analytical solution regardless of mesh resolution. \textcolor{black}{For
a  guide field of $B_{0}=10\gg B_{p}$, the error decreases as the
mesh is refined for sufficiently coarse meshes}, but it saturates
for finer meshes when the spatial discretization error becomes smaller
than that introduced by the tokamak-ordering approximation, which
for the magnetic field in Eq. \ref{magn_field} is of $\mathcal{O}(B_{0}^{-3})$
(as is evident from the figure by inspecting the error saturation
levels for different $B_{0}$). Increasing the guide field further
to $B_{0}=100\ggg B_{p}$ results in fourth-order-accurate convergence
for all meshes considered, although the scaling is expected to break
at some point for further refined meshes.

\subsection{Convergence study of the arbitrary-B formulation}

In contrast, the new arbitrary-B formulation is able to recover convergence
with the analytical solution for arbitrary guide-field strength $B_{0}$
and mesh resolution. For sufficiently large guide fields compared
to $B_{p}$, $B_{0}=10$, Fig. \ref{fig:tokamak_ordering_vs_arbB}-right
shows (for $\epsilon=10^{-2}$) that, unlike the tokamak-ordering
formulation, the new arbitrary-B formulation produces the expected
fourth-order asymptotic convergence rate for all meshes considered.

Convergence results for guide fields $B_{0}\lesssim B_{p}$ (i.e.,
more compressible magnetic field topologies) are more nuanced. This
is so because special magnetic-field topological surfaces such as
separatrices (which separate closed and open field lines) become more
relevant, impacting the accuracy of the simulation more. 
\begin{figure*}
\begin{centering}
\subfloat[\label{fig:arbB-anis2}$\epsilon=10^{-2}$]{\centering{}\includegraphics[width=0.48\columnwidth]{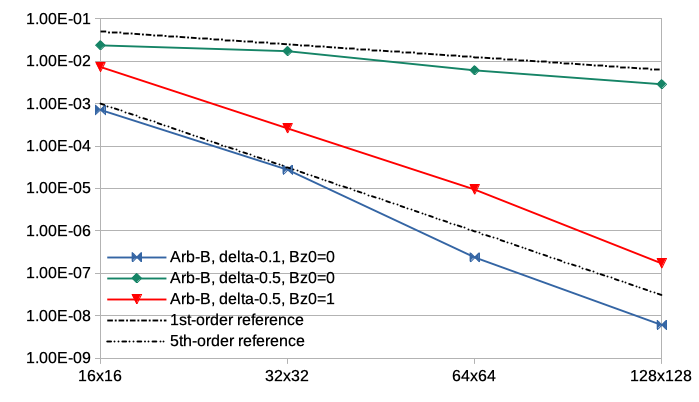}}\subfloat[\label{fig:arbB-anis4}$\epsilon=10^{-4}$]{\begin{centering}
\includegraphics[width=0.48\columnwidth]{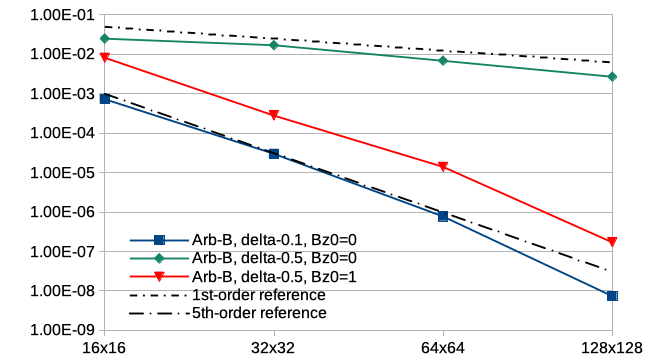}
\par\end{centering}
}
\par\end{centering}
\centering{}\subfloat[\label{fig:arbB-anis6}$\epsilon=10^{-6}$]{\begin{centering}
\includegraphics[width=0.48\columnwidth]{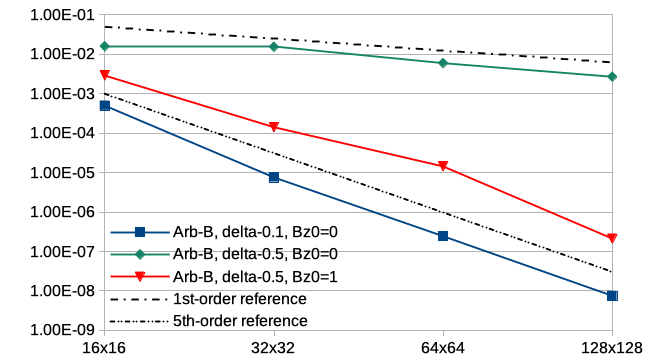}
\par\end{centering}
}\caption{Mesh convergence study with the arbitrary-B formulation for small
or zero magnetic guide fields.\label{fig:arb-B-conv-study}}
\end{figure*}
Fig. \ref{fig:arb-B-conv-study} depicts convergence results for $\epsilon=10^{-2},10^{-4},10^{-6}$
and various combinations of $B_{0}=0,1$ and $\delta=0.1,0.5$. The
first thing to notice is that error levels are largely independent
of the anisotropy ratio $\epsilon$, as expected from an AP formulation
\citep{chacon2014asymptotic,chacon-jpc-24-imp_xport}. Secondly, we
see order reduction of our numerical discretization to first-order
accuracy for the $B_{0}=0$, $\delta=0.5$ case, for which the verification
source is singular by construction due to the presence of magnetic
nulls. This is a consequence of the particular choice of verification
solution, and not a property of the method, which is otherwise able
to deal with it seamlessly in a convergent manner.

\begin{figure*}
\begin{centering}
\subfloat[\label{fig:arbB-anis4-d=00003D0.1-B0=00003D0} $\delta=0.1$ and $B_{0}=0$]{\centering{}\includegraphics[width=0.48\columnwidth]{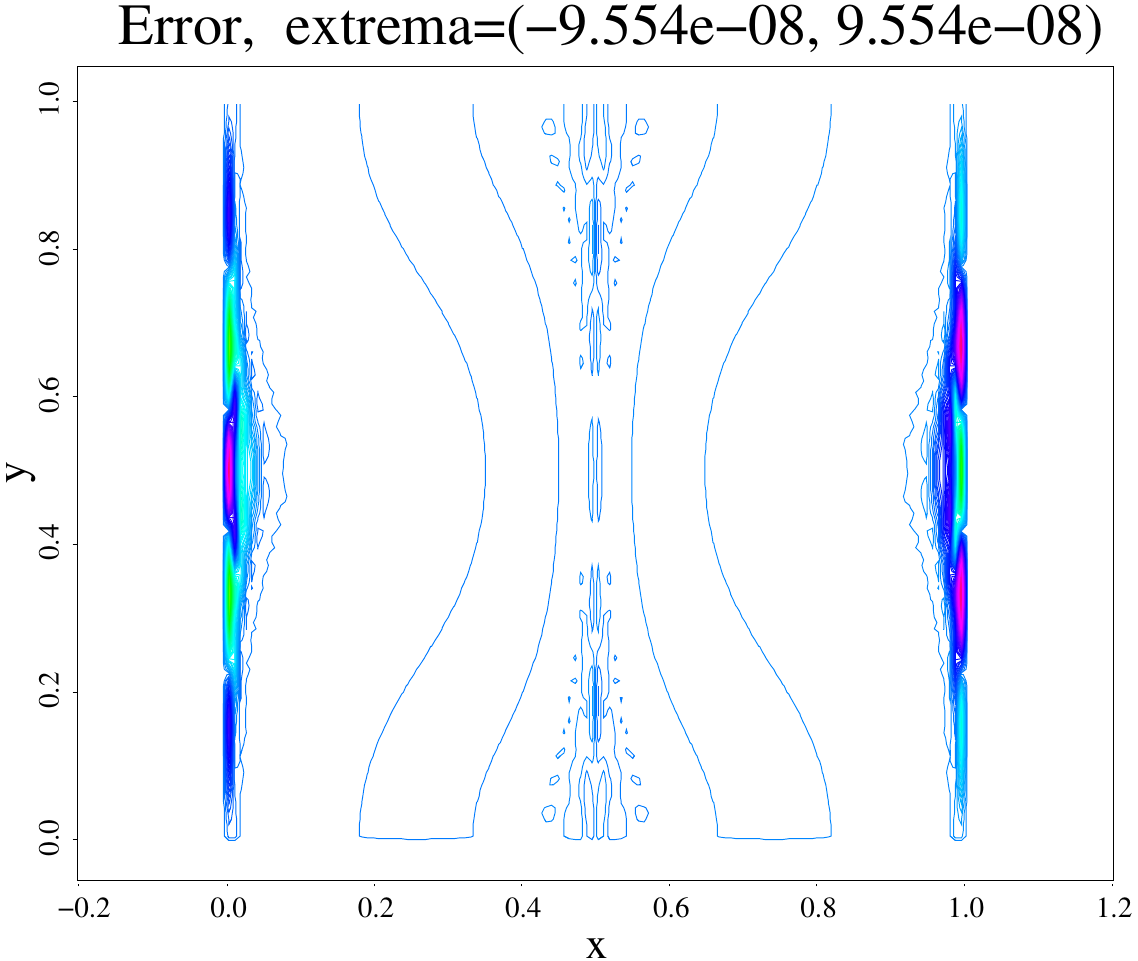}}\subfloat[\label{fig:arbB-anis4-d=00003D0.5-B0=00003D1} $\delta=0.5$ and $B_{0}=1$]{\begin{centering}
\includegraphics[width=0.48\columnwidth]{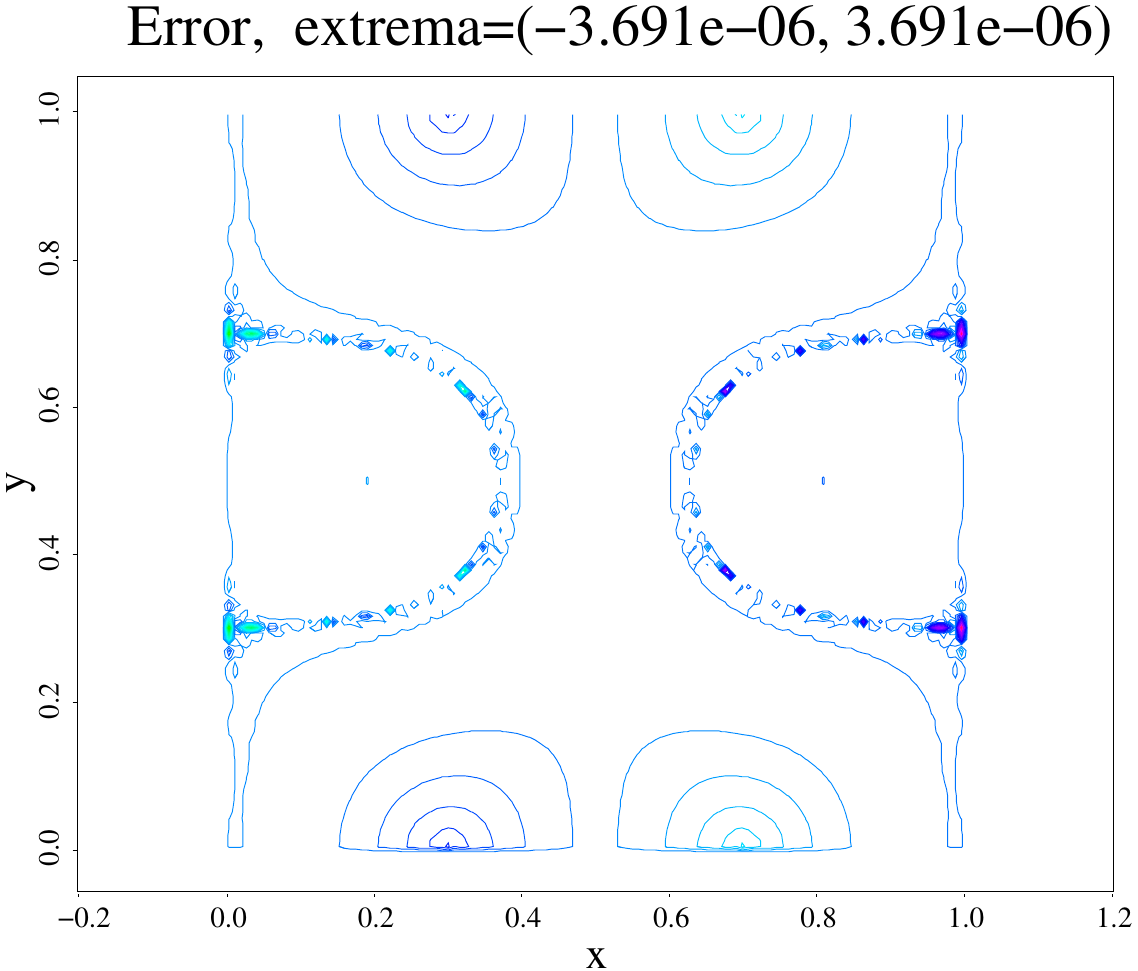}
\par\end{centering}
}
\par\end{centering}
\centering{}\subfloat[\label{fig:arbB-anis4-d=00003D0.1-B0=00003D10} $\delta=0.1$ and
$B_{0}=10$]{\centering{}\includegraphics[width=0.48\columnwidth]{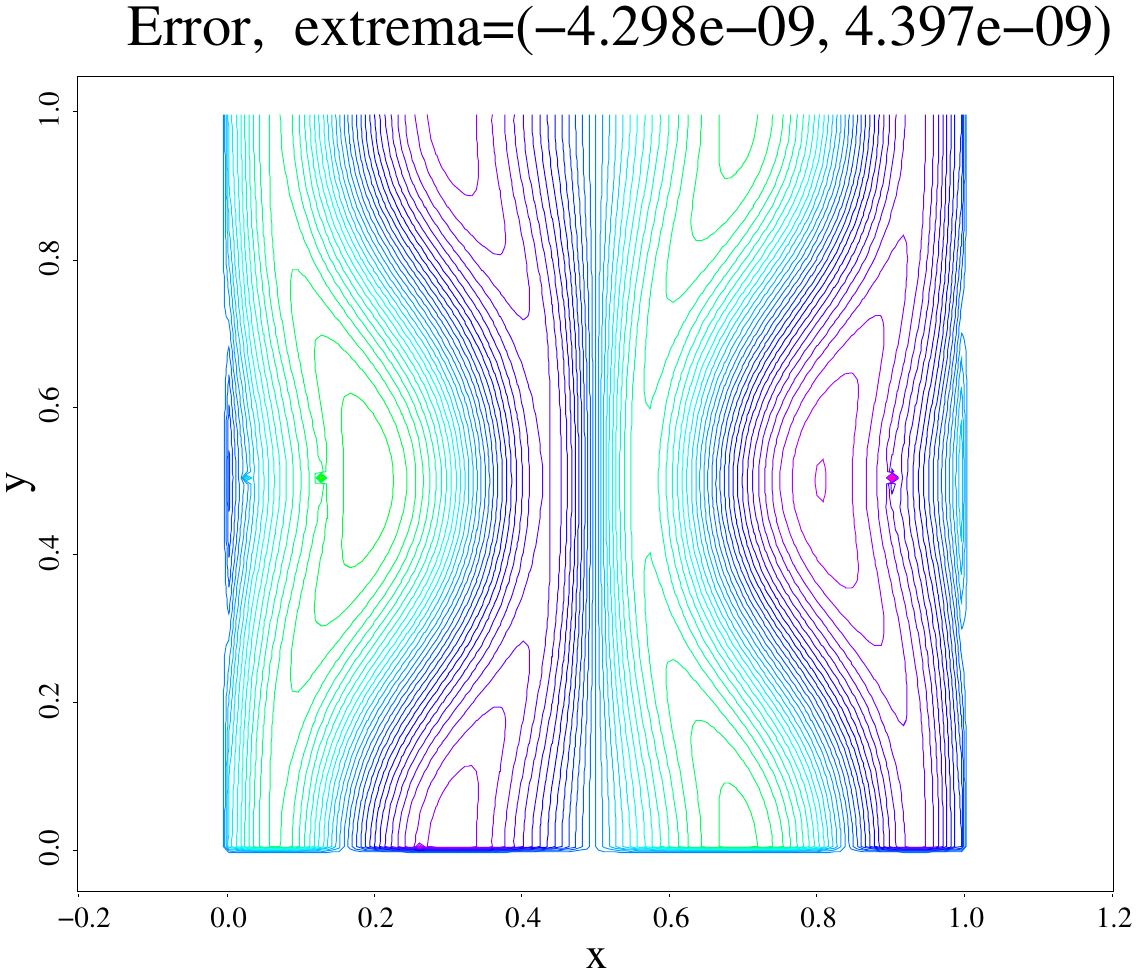}}\subfloat[\label{fig:arbB-anis4-d=00003D0.5-B0=00003D10} $\delta=0.5$ and
$B_{0}=10$]{\begin{centering}
\includegraphics[width=0.48\columnwidth]{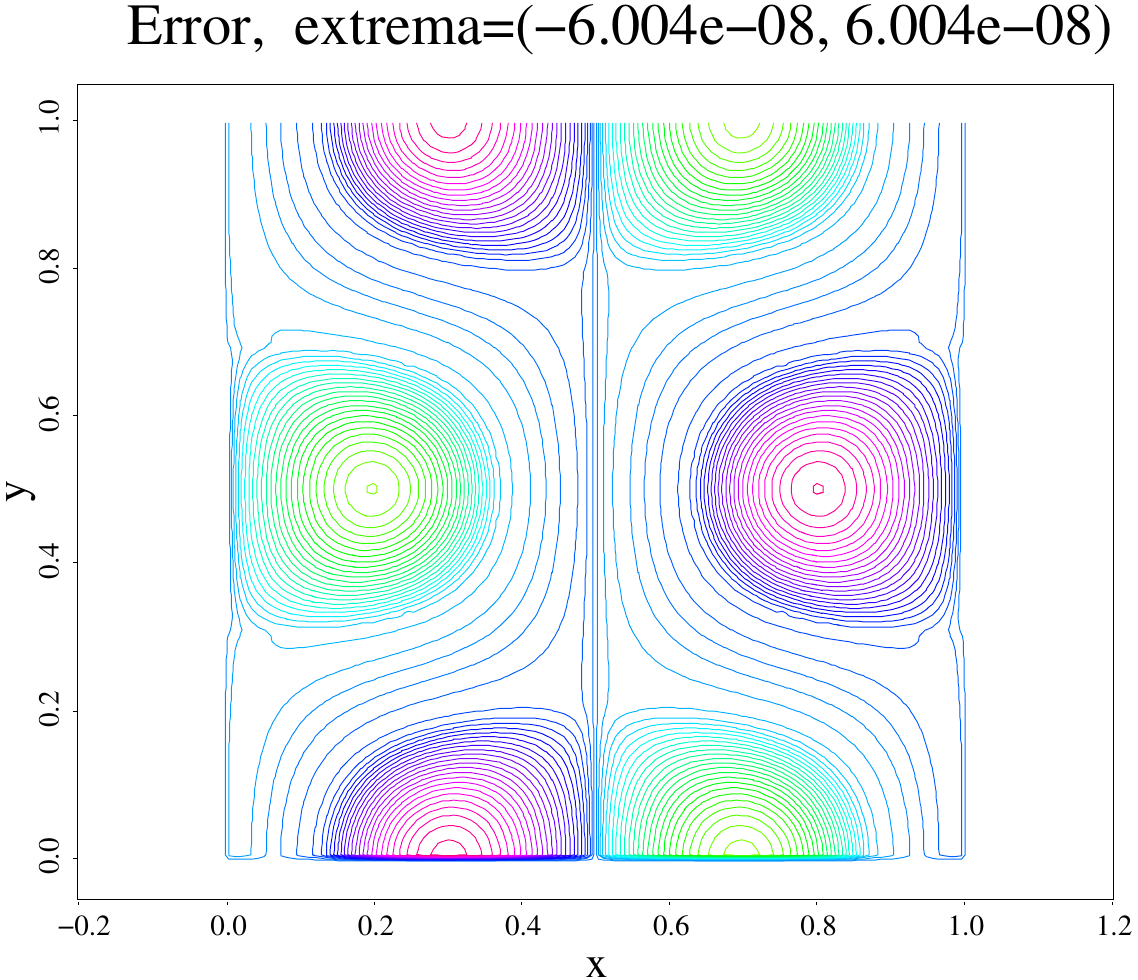}
\par\end{centering}
}\caption{Steady-state errors vs the analytical solution for $\epsilon=10^{-4}$
on the $128\times128$ mesh for various choices of $\delta$ and $B_{0}$.\label{fig:arb-B-error}}
\end{figure*}

In the cases where the verification source is regular but still with
$B_{0}\lesssim B_{p}$ (e.g., $\delta=0.1$, $B_{0}=0$ or $\delta=0.5$,
$B_{0}=1$), we actually find super-convergence with mesh refinement,
of fifth-order instead of the expected fourth-order. Further inspection
of the numerical errors for these cases (see Figs. \ref{fig:arbB-anis4-d=00003D0.1-B0=00003D0}
and \ref{fig:arbB-anis4-d=00003D0.5-B0=00003D1} for $\epsilon=10^{-4}$)
reveals that the errors for $B_{0}\lesssim B_{p}$ concentrate along
one-dimensional manifolds, e.g. along the boundary for $\delta=0.1$
and around magnetic-island separatrices for $\delta=0.5$, which explains
the super-convergence (since the $\ell_{2}$-norm accounts for volume
elements in the computation of the errors, Eq. \ref{eq:l2-norm}).
These 1D manifolds are locations where boundary layers form, and where
the numerical error is expected (and found) to be largest \citep{chacon2014asymptotic,chacon-jpc-24-imp_xport}.
The 1D error manifolds disappear as one increases the guide field
sufficiently, as evidenced by Figs. \ref{fig:arbB-anis4-d=00003D0.1-B0=00003D10}
and \ref{fig:arbB-anis4-d=00003D0.5-B0=00003D10}, which are obtained
for $\epsilon=10^{-4}$ and $B_{0}=10\gg B_{p}$.

\section{Discussion and conclusions}

\label{sec:conclusions}We have proposed an extension of the recently
proposed semi-Lagrangian formulation \citep{chacon2014asymptotic,chacon-jpc-24-imp_xport}
for the highly anisotropic heat transport equation to deal with arbitrary
magnetic field configurations by lifting the so-called tokamak-ordering
approximation, which required $\nabla\cdot(\mathbf{B}/B)\sim0$. The
new approach still leverages analytical Green's functions for field-line
integrals (which is critical to make the method tractable) by employing
a suitable change of variables and introducing \textcolor{black}{a suitable
positive-definite field $\beta(\mathbf{x})$}. A prescription for
$\beta(\mathbf{x})$ is proposed that renders the discrete formulation
consistent, AP, and able to deal with magnetic nulls seamlessly. The
method has been verified with a manufactured solution for arbitrary
guide fields and anisotropy ratios, and has been demonstrated to achieve
at least fourth-order design spatial accuracy in all cases without
magnetic field nulls. Magnetic field nulls render the verification
solution source singular, resulting in order reduction to first order,
but does not break the method. This order reduction is a consequence
of the particular choice of verification solution rather than a feature
of the method. Future extensions of this work will consider nonlinear
heat transport coefficients and finite magnetic field lines, e.g.
when magnetic field lines intersect the domain boundary. 

\section*{Acknowledgments}

This work was supported by Triad National Security, LLC under contract
89233218CNA000001 and DOE Office of Applied Scientific Computing Research
(ASCR). The research used computing resources provided by the Los
Alamos National Laboratory Institutional Computing Program. LC acknowledges
useful conversations with C. Hauck, D. del-Castillo-Negrete, and O.
Koshkarov. GDG participated in this research while visiting LANL in
Spring 2019.

\appendix

\section{Manufactured solution source}

\label{sec:Manufactured-source}

The manufactured solution source is found from:
\begin{equation}
S=-(1-\epsilon)\nabla_{\parallel}^{2}\tilde{T}-\epsilon\nabla^{2}\tilde{T}-\nabla^{2}T_{\mathcal{N}}(\psi).
\end{equation}
Using (per the main text) $T_{0}(\psi)=x+\delta\sin(2\pi x)\cos(2\pi y)$
and $\tilde{T}=\cos(4\pi y)\sin(4\psi x)$, we find:

\begin{eqnarray*}
\nabla^{2}\tilde{T} & = & -32\pi^{2}\sin(4\pi x)\cos(4\pi y)\\
\nabla^{2}T_{\mathcal{N}}(\psi) & = & -8\pi^{2}\delta\sin(2\pi x)\cos(2\pi y)\\
\nabla_{\parallel}^{2}\tilde{T} & = & 4\pi\Biggl\{-\frac{4\pi^{2}}{B^{2}}\biggl(\Bigl(B_{x}^{2}+B_{y}^{2}\Bigr)\cos(4\pi y)\sin(4\pi x)\\
 &  & +2B_{x}B_{y}\cos(4\pi x)\sin(4\pi y)\biggr)\\
 &  & +\;\cos(4\pi x)\cos(4\pi y)\left(\partial_{y}\left(\frac{B_{x}B_{y}}{B^{2}}\right)+\partial_{x}\left(\frac{B_{x}^{2}}{B^{2}}\right)\right)\\
 &  & -\;\sin(4\pi x)\sin(4\pi y)\left(\partial_{x}\left(\frac{B_{x}B_{y}}{B^{2}}\right)+\partial_{y}\left(\frac{B_{y}^{2}}{B^{2}}\right)\right)\Biggr\}.
\end{eqnarray*}
The derivative terms read:

\begin{eqnarray*}
\partial_{x}\left(\frac{B_{x}^{2}}{B^{2}}\right) & = & 2\mathcal{A}B_{x}\biggl\{\cos(2\pi x)\sin(2\pi y)\left(1-\frac{B_{x}}{B^{2}}\right)\\
 & \; & +\sin(2\pi x)\cos(2\pi y)\frac{B_{x}B_{y}}{B^{2}}\biggr\}\\
\partial_{x}\left(\frac{B_{x}B_{y}}{B^{2}}\right) & = & \mathcal{A}\Biggl\{ B_{y}\cos(2\pi x)\sin(2\pi y)\left(1-2\frac{B_{x}^{2}}{B^{2}}\right)\\
 & \; & -B_{x}\cos(2\pi y)\sin(2\pi x)\left(1-2\frac{B_{y}^{2}}{B^{2}}\right)\Biggr\}\\
\partial_{y}\left(\frac{B_{y}^{2}}{B^{2}}\right) & = & 2\mathcal{A}B_{y}\biggl\{\sin(2\pi y)\cos(2\pi x)\left(\frac{B_{y}^{2}}{B^{2}}-1\right)\\
 &  & -\cos(2\pi y)\sin(2\pi x)\frac{B_{x}B_{y}}{B^{2}}\biggr\}\\
\partial_{y}\left(\frac{B_{x}B_{y}}{B^{2}}\right) & = & \mathcal{A}\Biggl\{ B_{y}\cos(2\pi y)\sin(2\pi x)\left(1-\frac{2B_{x}^{2}}{B^{2}}\right)\\
 &  & -B_{x}\cos(2\pi x)\sin(2\pi y)\left(1-\frac{2B_{y}^{2}}{B^{2}}\right)\Biggr\},
\end{eqnarray*}
with $\mathcal{A}=\frac{4\pi^{2}\delta}{B^{2}}$ . Note that the source
is singular for magnetic nulls, $B=0$.

\section{\textcolor{black}{Error of the operator-split arbitrary-B Lagrangian
formulation}}

\label{app:os-error}

\textcolor{black}{An operator-split algorithm can be formulated from
the formal result in Eq. \ref{eq:formal_solution_general_B} as suggested
in Ref. \citep{chacon2014asymptotic} by considering first an update
of the slow dynamics (which exactly respects the null space \citep{chacon2014asymptotic}):
\begin{equation}
\frac{T^{*}-T^{n}}{\Delta t}=\nabla_{\perp}^{2}T^{*}+S,\label{eq:split_slow_general_B-2}
\end{equation}
followed by a Lagrangian step in which the new-time temperature in
the source $S^{*}$ (Eq. \ref{eq:formal_source_beta}) is approximated
by $T^{*}$. It follows that $S^{*}=\beta\frac{T^{*}-T^{n}}{\Delta t},$
and hence the Lagrangian splitting stage simply reads:
\begin{eqnarray}
T_{OS}^{n+1}(\mathbf{x}) & = & \mathcal{G}_{\lambda}\left(T^{n};\mathbf{x},\frac{\Delta t}{\beta\epsilon}\right)+\mathcal{P}_{\lambda}\left(T^{*}-T^{n};\mathbf{x},\frac{\Delta t}{\beta\epsilon}\right)+\mathcal{E}_{split},\label{eq:split-fast-general_B-2}
\end{eqnarray}
which is otherwise identical to the one proposed in Ref. \citep{chacon2014asymptotic}
except for the parameter $\beta$, and the fact that the Lagrangian
integrals are performed in the variable $\lambda$.}

\textcolor{black}{However, as formulated, this operator-split algorithm
is }\textcolor{black}{\emph{not}}\textcolor{black}{{} convergent with
either small $\Delta t$ or small $\epsilon$. This can be readily
shown by computing the operator-splitting local-truncation error (LTE)
by subtracting Eqs. \ref{eq:split-fast-general_B-2} and \ref{eq:formal_solution_general_B},
which can be written after some manipulation as:
\[
\mathcal{E}_{split}=T^{n+1}-T_{OS}^{n+1}\approx\mathcal{P}_{\lambda}\left((T^{n+1}-T^{*})-\frac{T^{n+1}-T^{n}-\Delta t(\nabla_{\perp}^{2}T^{n+1}+S)}{\beta B^{2}};\mathbf{x},\frac{\Delta t}{\beta\epsilon}\right).
\]
Using the BDF1 semi-discretization of the original transport PDE,
we can write this result as:
\[
\mathcal{E}_{split}=T^{n+1}-T_{OS}^{n+1}\approx\mathcal{P}_{\lambda}\left((T^{n+1}-T^{*})-\frac{\Delta t}{\epsilon\beta}\partial_{\lambda}^{2}T^{n+1};\mathbf{x},\frac{\Delta t}{\beta\epsilon}\right).
\]
For the first splitting-error contribution, $\mathcal{P}_{\lambda}\left(T^{n+1}-T^{*};\mathbf{x},\frac{\Delta t}{\beta\epsilon}\right)$,
a Fourier analysis yields (following Ref. \citep{chacon2014asymptotic}):
\[
\hat{\mathcal{E}}_{split,I}^{k_{\parallel}}\sim\hat{T}_{k_{\parallel}}\Delta t\min\left[1,\frac{k_{\parallel}^{2}\Delta t}{\epsilon}\right],
\]
which is convergent in all regimes of $\Delta t/\epsilon$. However,
for the second contribution, we have:
\[
\mathcal{E}_{split,II}=\frac{\Delta t}{\epsilon\beta}\mathcal{P}_{\lambda}\left(\partial_{\lambda}^{2}T^{n+1};\mathbf{x},\frac{\Delta t}{\beta\epsilon}\right)=\mathcal{G}_{\lambda}\left(T^{n+1};\mathbf{x},\frac{\Delta t}{\beta\epsilon}\right),
\]
which follows by integration by parts of the kernel integral and using
that $\partial_{s}^{2}\mathcal{U}(s,t)=\frac{1}{t}G(s,t)$. Therefore,
we find:
\[
\mathcal{E}_{split,II}=\mathcal{G}_{\lambda}\left(T^{n+1};\mathbf{x},\frac{\Delta t}{\beta\epsilon}\right)\rightarrow\begin{cases}
T^{n+1}(\mathbf{x}), & \Delta t\ll\beta\epsilon\\
\left\langle T^{n+1}\right\rangle , & \Delta t\gg\beta\epsilon
\end{cases},
\]
which does not vanish with either arbitrarily small $\Delta t$ or
$\epsilon$, and is therefore not convergent.}

\pagebreak{}

\bibliographystyle{ieeetr}
\bibliography{ref}

\end{document}